\documentclass[11pt]{article}                                                   
\renewcommand{\thefootnote}{\fnsymbol{footnote}}
\renewcommand{\thanks}[1]{\footnote{#1}} 
\newcommand{\starttext}{
\setcounter{footnote}{0}
\renewcommand{\thefootnote}{\arabic{footnote}}}

\newcommand{\be}{\begin{equation}}
\newcommand{\bea}{\begin{eqnarray}}
\newcommand{\eea}{\end{eqnarray}}
\newcommand{\ee}{\end{equation}}

\def\ba{\begin{eqnarray}}
\def\ea{\end{eqnarray}}

\def\G{{\cal G}}
\def\X{{\cal X}}

\def\p{\partial}

\def\Chow{{\rm Chow}}

\def\o{\omega}

\def\l{\lambda}
\def\m{\mu}
\def\n{\nu}
\def\o{\omega}

\def\t{\theta}

\def\Z{{\bf Z}}
\def\Q{{\bf Q}}
\def\R{{\bf R}}
\def\C{{\bf C}}
\def\P{{\bf P}}
\def\cX{{\cal X}}
\def\cY{{\cal Y}}
\def\cK{{\cal K}}
\def\cL{{\cal L}}
\def\cM{{\cal M}}
\def\cZ{{\cal Z}}
\def\cS{{\cal S}}

\def\Fut{{\rm Fut}}
\def\AY{{\rm AY}}
\def\si{\sigma}

\def\i{\infty}
\def\I{\int}
\def\p{\prod}
\def\s{\sum}

\def\sub{\subseteq}
\def\ra{\rightarrow}
\def\ddb{\partial\bar\partial}

\def\cF{{\cal F}}
\def\cA{{\cal A}}

\def\G{\Gamma}

\def\cN{{\cal N}}

\def\v{\medskip}

\def\[{{\bf [}}
\def\]{{\bf ]}}

\def\pl{\partial}
\def\bo{{\bf 1}}
\def\cO{{\cal O}}

\begin{document}
\starttext
\setcounter{footnote}{0}

\title{The Futaki Invariant and the  Mabuchi Energy of a Complete Intersection}

\author{D.H. Phong and Jacob Sturm
}

\date{}

\maketitle
\footnotetext[1]{Research supported in part by the National
Science Foundation under grants
DMS-02-45371 and DMS-01-00410.}

\setcounter{equation}{0}
\setcounter{footnote}{0}
\textwidth=5in                                                                  
\textheight=7.5in

\section{Introduction.}

{
Let $M$ be a compact complex K\"ahler manifold. If $c_1(M)=0$ or
if $c_1(M) < 0$, then it is known by the work of Yau \cite{Yau78} and
Yau, Aubin \cite{Yau78}, \cite{Aubin} that $M$ has a K\"ahler-Einstein
metric. If $c_1(M)>0$, then there are obstructions to the existence
of such a metric, and here the guiding conjecture is that formulated by
Yau in \cite{Yau93}, which says that $M$ has a K\"ahler-Einstein metric
if and only if $M$ is stable in the sense of geometric invariant theory.

\medskip

An important obstruction to the existence of  K\"ahler-Einstein metric
is the invariant of Futaki \cite{Fut 83}, which is a map 
$F: \eta(M)\ra \C$ with the following properties: $F$ is
a Lie algebra character on the space $\eta(M)$ of holomorphic vector fields, which depends only on the cohomology class
$[\o]\in H^2(M)$. The vanishing of $F$ is a necessary
condition for the existence of a K\"ahler-Einstein metric on $M$.
However, it is not a sufficient condition: in \cite{T97},  Tian 
gives an example of a manifold with
$\eta(M)=0$ (so that the Futaki invariant vanishes trivially) with the property that $M$ has no K\"ahler-Einstein metric.

\medskip

In \cite{DT}, a generalized Futaki invariant 
is defined for Q-Fano varieties, which are certain singular
varieties arising naturally
as degenerations of smooth Fano manifolds. In
\cite{T97}, a manifold $M$ is defined to be K-stable if the Futaki invariant of every
non-trivial Q-Fano  degeneration of $M$ has positive real
part. It is proved there that if $M$ has a
K\"ahler-Einstein metric, then
$M$ is  K-stable, and the converse is conjectured
to be true as well.

\medskip

In general, it is rather difficult to check the K-stability of a given
manifold. Recent progress has been made by Lu \cite{Lu2} in the case where $M$ is a hypersurface: Using a delicate analytic argument, he  provides,
in the
case where $M$ is a hypersurface, an explicit formula for the Futaki invariant
of a degeneration of $M$. This represents a non-linear generalization
in the special case of hypersurfaces of
Lu's earlier formula \cite{Lu} for the Futaki invariant of a normal complete intersection.  

\medskip

In this paper we shall study the Futaki invariant and the Mabuchi K-energy using the Deligne pairing technique. The Deligne
pairing technique
was introduced by Zhang \cite{Z} and developed 
in \cite{PS1} and \cite{PS2}. Our main result
is a formula (Theorem 6) for the Mabuchi energy on the
orbits of $SL(N+1)$ which provides a non-linear version of Lu's Futaki invariant formula for a complete intersection. 
A similar formula is also established for the Aubin-Yau
functional. We also give a simple characterization of the generalized Futaki
invariant and a new construction of the generalized Futaki character as well. A basic idea in our approach is to
associate to each smooth projective variety $M$ the one-dimensional vector space
\be
{\cal F}=\langle K^{-1},\cdots,K^{-1}\rangle
\ee
Here $K$ is the canonical bundle of $M$, $\langle\cdots\rangle$
is the Deligne intersection pairing, and $K^{-1}$ occurs
$n+1$ times in the pairing. The point is that this construction
is completely canonical, and hence the action of $\eta(M)$ on $M$
lifts to a group action on ${\cal F}$, which turns out to
be precisely the Futaki character. A slightly modified version
works in the Q-Fano case as well.

\medskip

A more detailed description of the results 
and organization of the paper is as follows. 
Section 2 is devoted to establishing the
simple characterization of the Futaki invariant: $F(X)$ is 
the eigenvalue of the infinitesimal action of $X$ on
$\Chow(M)$, the Chow point of
$M$. Moreover, the Futaki character $\hat F(X)$ is the eigenvalue of the action
of Aut$(M)$ on $\Chow(M)$
(see Theorem 1 and its corollary, stated in \S 2.1 and proved in
\S 2.4).
The background needed on Deligne pairings is provided in
\S 2.2. Energy functionals in K\"ahler geometry
originate from many sources. A unifying theme is that of
Bott-Chern secondary classes (see e.g. \cite{Donaldson},
\cite{FMS}, \cite{T00} and references therein).
In this paper, we rely on another unifying theme
discussed in \S 2.3,
namely energy functionals as variations of metrics in suitable Deligne pairings 
(Theorem~2), e.g.
\be
{\cal M}_{he^{-\phi}}={\cal M}_h\otimes O(\nu_{\o}(\phi))
\ee  
Next, in Section 3, we show how one can use this
characterization to give a new proof of Lu's formula
\cite{Lu} on the generalized Futaki invariant of a complete intersection (Theorem 3). Section 4
begins with Theorem 4,
which provides a general ``adjunction formula with metrics". This
theorem is closely related to a corresponding result of Lu, although the
formulation in \cite{Lu} is rather different. Finally, we apply the adjunction
formula with metrics in Theorems 5 and 6 to prove 
the desired formulas for the Aubin-Yau and the
Mabuchi energy functionals for complete intersections.
As stated earlier, these can be viewed as a non-linear version 
of Lu's formula for the Futaki
invariant of complete intersections.
They can also be viewed as 
a generalization to the case of complete
intersections of formulas for hypersurfaces
obtained earlier in
\cite{T94} and \cite{PS1} for the Mabuchi functional,
and in \cite{Z} and \cite{PS1} for the Aubin-Yau
functional: in particular, we show that the K-energy of a complete intersection may be
expressed as a degenerate norm on the space of defining polynomials. 
   }

\section{Lifting the Futaki character}
\setcounter{equation}{0}

\subsection{Statement.}

Let $M$ be an irreducible normal projective variety and $M_{reg}\sub M $
the open subset consisting of all smooth points. We
say
$M$ is Q-Fano if there exists $k>0$ and a very ample line bundle
$L$ on $M$ such that $L|_{M_{reg}}= K_{M_{reg}}^{-k}$. A K\"ahler
form $\o$ on $M_{reg}$ is {\it admissible} if there exists
$\phi_{L^m}:M\hookrightarrow \P^N$, an embedding of $M$ defined
by some power of $L$, and a K\"ahler form $\tilde\o$ on $\P^N$
representing $c_1(\P^N)$, such that 
\be \label{omega}
\o\ = \ {1\over km}\phi^*_{L^m}(\tilde\o)
\ee
If $\tilde h$ is a metric on $O(1)$ such that $
\tilde\o=Ric(\tilde h)
\equiv -{\sqrt{-1}\over 2\pi}
\ddb \tilde h$,
then we define a metric $h$ on $L$ as follows:
\be \label{h}
h \ = \ {1\over m}\phi^*_{L^m}(\tilde h)
\ee

A vector field $X$ on $M_{reg}$ is admissible if there exists
an embedding as above, and a holomorphic vector field $\xi$
on $\P^N$ such that $\xi$ is tangent to $M_{reg}$ with 
$\xi|_{M_{reg}}=X$. The space of all admissible vector fields
will be denoted $\eta(M)$. It is the Lie algebra of the
algebraic group $G=Aut(M)$ consisting of all invertible
algebraic maps from $M$ to itself.

\medskip

Suppose $\o$ is an admissible K\"ahler form on $M_{reg}$ and
$X\in\eta(M)$. Define
\be \label{futaki}
F(X)\ = \ (n+1)\I_{M_{reg}} X(f)\o^n
\ee
where $f$ is a smooth function on $M_{reg}$ chosen so that
$$ Ric(\o)-\o\ = \ {\sqrt{-1}\over 2\pi} \pl\bar\pl f
$$
If $M$ is smooth, then Futaki \cite{Fut 83}
has shown that (\ref{futaki}) is well defined,
independent of the choice of~$\o$, and $F:\eta(M)\ra \C$ is
a Lie algebra character. For arbitrary Q-Fano $M$, this was proved, using
resolution of singularities,
by  Ding-Tian \cite{DT}. When $M$ is
smooth, it is known, through the work of Futaki, Mabuchi and Morita (see
\cite{Fut 87 Inv.},
\cite{Mab K energy maps},
 \cite{Fut-Morita}) that $F$ has a lift to a
group character
$\hat F: G\ra\C^\times$. They show, using the theory of Chern-Simons
invariants, that
$\hat F$ may be constructed as the fiber integral of the Godbillon-Vey
class of a certain locally trivial $M$ bundle over the classifying space
$BG$. Yotov \cite{Y2} was able to construct $\hat F$ in the general
case of Q-Fano varieties by 
combining the techniques of the previous authors with the
Edidin-Graham \cite{EG} theory of equivariant Chow cohomology groups.

\medskip

Our first theorem gives an alternate
construction of $\hat F$ in terms of the Deligne pairing
(for background on Deligne pairings, see \S 2.3): Let
$M$ be a Q-Fano variety and let
$L$ be the extension of $K^{-k}$ from $M_{reg}$ to $M$. Then
the action of
$G$ on
$M$ lifts canonically to an action on $K^{-1}$ and thus $G$ acts on the
pairing
$\cF = \langle L,L,...,L \rangle$, the Deligne pairing of 
$L$ with itself $n+1$ times. Since $\cF$ is a one
dimensional vector space,  the action of $G$ on $\cF$ defines a group
character 
$$\hat F: G \ra Aut (\cF)=\C^\times $$ Thus, if
$\ell_0,...,\ell_n$ are rational sections of $L$ with
$div(\ell_0)\cap\cdots\cap div(\ell_n)=\emptyset$, 
and $\sigma\in G$, then
$\hat F(\sigma)= \langle \si_*\ell_0,...,\si_*\ell_n\rangle/
 \langle \ell_0,...,\ell_n\rangle  \in \C^\times.$
\v

{\bf Theorem 1.}{\it \ The character $\hat F: G \ra Aut
(\cF)=\C^\times $ is a lifting of the Futaki character. More precisely,
if $X\in\ h(M)$  and if  $X_\R= X+\bar X$ is the real part of $X$, then
$$ \hat F(\exp(tX_\R))\ = \ \exp(t\,k^{n+1}F(X))
$$
}As a corollary, we find that $F(X)$ is the eigenvalue of the
Chow point of $M$: More, precisely, there is a natural action
of $G$ on $H^0(L^k)$, inducing an action on $Gr(N-n-1,\P^N)$ and
$H^0(Gr,O(d))$, where $d$ is the degree of $M\sub\P^N$, and
$Gr=Gr(N-n-1,\P^N)$ is the Grassmannian manifold consisting of
all $N-n-1$ planes in $\P^N$. Let $\Chow_k(M)\in H^0(Gr,O(d))$ be
a representative of the Chow point. Since the Chow point is
unique up to scalar multiplication, we have, for each $\sigma\in G$,
$$ \sigma\cdot \Chow_k(M)\ \equiv \Chow_k(M)\circ\sigma^{-1} \ = \ 
C(\sigma)
\Chow_k(M)
$$
for some $C(\sigma)\in\C^\times$. The corollary states that $C(\sigma)=
\hat F(\sigma)^{k^{n+1}}$:
\v

{\bf Corollary 1.} {\it The Futaki group character $\hat F(\sigma)$ is characterized by
$$ \sigma\cdot \Chow_k(M)\ = \ \hat F(\sigma)^{k^{n+1}}  \Chow_k(M)
$$
and the Futaki
invariant $F(X)$ by
$$ X\cdot Chow_k(M)\ = \ k^{n+1}\,F(X)\, \Chow_k(M)
$$
for all $\sigma\in G\sub GL(N+1,\C)$, and for all $X\in Lie(G)$.
}
\v

\subsection{The Deligne pairing}

We recall some of the basic definitions and properties in \cite{De} and
\cite{Z}
\footnote{At the December 2002 Complex Geometry conference in Tokyo, Professors T. Mabuchi and L. Weng informed us that
they have also been aware for some time of potential applications of Deligne pairings to K\"ahler geometry.}: Let $\pi:\cX \ra S$ be a flat projective morphism of integral
schemes of relative dimension~$n$. Thus for every $s \in S$,
the fiber $\cX_s$ is a projective variety in $\P^N$ of dimension
$n$.
Let $\cL_0, \cL_1,...,\cL_n$ be line bundles on $\cX$. The
Deligne pairing is a line bundle on $S$, denoted
$\langle\cL_0,\cL_1,...,\cL_n\rangle(\cX/S)$,
and defined as follows:
Let $U \sub S$ be a small open set and let
$l_i$ be a rational section of $\cL_i$ over $\pi^{-1}U$.
Thus $l_i: \pi^{-1}(U) \ra \cL_i$ is a rational function
with the property $p_i\circ l_i = \bo_{\pi^{-1}U}$, where
$p_i:\cL_i \ra \cX$ is the projection map of $\cL_i$.
Assume that the $l_i$ are chosen in ``general
position": This means $\cap_i\  div(l_i)=\emptyset $
and for each $s$ and $i$ with
$s\in U$ and $0\leq i \leq N$, the fiber
 $\cX_s$ is not
contained in $div(l_i)$.
 Then for every $k$, the map
$\big(\cap_{i\not= k} \ div(l_i)\big) \ra S$ is finite:
For
every $s$, 
$\big(\cap_{i\not= k} \ div(l_i)\big)\cap \cX_s \sub \cX_s$ is
a zero cycle $\s n(s)P(s) $. This means  that the $n(s)$ are
integers and the points $P(s)$ are a finite set of points in $\cX_s$.
\v

Now we define $\langle\cL_0,\cL_1,...,\cL_n\rangle(\cX/S)$. Over a
small $U \sub S$, this line bundle is trivial and
generated by the symbol $\langle l_0,...,l_n\rangle$ where the
$l_i$ are chosen to be in general position. If
$l_i'$ is another set of rational sections in
general position, then
$\langle l_0',...,l_n'\rangle = \psi(s)\langle l_0,...,l_n\rangle$ for some
nowhere vanishing function $\psi$ on $U$ which we must specify.
We do this one section at a time: Assume that $l_i = l_i'$
for all $i \not=k$. Assume as well that the
rational function $f_k = l_k'/l_k$ is well defined
and non-zero on
$\big(\cap_{i\not= k} \ div(l_i)\big)$. Then
$ \psi(s) \ = \ \p f(P(s))^{n(s)}$.
The Deligne pairing is obviously multilinear in
the bundles $\cL_0,\cdots,\cL_n$.
Applied to a family $\cX$ consisting of a single irreducible
normal variety (i.e. $S$ consists of a single point), it produces a one-dimensional vector space.
\v

Let $\pi: \cX \ra S$ and $\cL_0, ..., \cL_n$ as above.
Let $l$ be a rational section of $\cL_n$. Assume
all components of $div(l)$ are flat over $S$.
Then we have  the following
{\it  induction formula: }
\be \label{1.1.5} \langle\cL_0,...,\cL_n\rangle(\cX/S) \ = \
\langle \cL_0,...,\cL_{n-1}\rangle(div(l)/S)
\ee
where $\langle \cL_0,...,\cL_{n-1}\rangle(Z/S)
=\prod_i\langle \cL_0,...,\cL_{n-1}\rangle(Y_i/S)^{n_i}$
if $Z=\sum_in_iY_i$ is a cycle on $S$.

\v

Assume now that $\cX,S$ are defined over $\C$, and
that $\cL_i$ is endowed with a smooth hermitian metric
(that is, for any holomorphic map from a smooth variety $Y$
to $\cX$, the pull back metric on $f^*\cL$ is smooth on $Y$).
We now define a hermitian metric on
$\langle\cL_0,...,\cL_n\rangle(\cX/S)$: Let 

$$c_1'(\cL_i)=-{\sqrt{-1}\over
2\pi}\ddb \log ||l||^2$$ be the normalized  curvature of $\cL_i$
where $l$ is a local invertible
section of $\cL_i$. Then $c_1'(\cL_i)$ is a $(1,1)$ form
 on $\cX_{reg}$.
The metric on $\langle\cL_0,...,\cL_n\rangle(\cX/S)$ is defined
by induction:
When $n=0$,
$\langle\cL\rangle_s \ =  \otimes_{p \in \pi^{-1}(s)  } \ \cL_p $
so we define
\be\label{2.2}
|| \langle l_0\rangle||_s \ =
 \ \p_{p \in \pi^{-1}(s)  }
||l_0(p)||
\ee
In general, we define
$$\label{2.3}
\log ||\langle l_0,..., l_n\rangle|| \ = \
\log||\langle l_0,...,l_{n-1}\rangle (div(l_n)/S)|| \ + \
\I_{\cX/S} \log ||l_n|| \Lambda_{i=0}^{n-1}c'_1(\cL_i)
$$
where the integral is the fiber integral over $S$,
so both sides  are functions on the base manifold $S$.

\v

If we combine the induction formula
with the definition of the metric,
we immediately get the following isometry:
$$
\langle \cL_0,...,\cL_n\rangle(\cX/S) \ = 
$$
\be\label{1.2.1}\
\langle\cL_0,...,\cL_{n-1}\rangle(div(l)/S)\otimes \cO
\left(-\I_{\cX/S} \log ||l_n||
\Lambda_{i=0}^{n-1}c'_1(\cL_i)\right)
\ee
where $\cO(f)$ denotes the trivial line bundle with
metric $||1|| = \exp(-f)$.
In particular,
\be \label{small change}
\langle\cL_0,...,\cL_{n-1},\cL_n\otimes\cO(\phi)\rangle(\cX/S) \ = \
\langle\cL_0,...,\cL_n\rangle(\cX/S)\otimes \cO(E)
\ee
where
\be\label{2.6}  E \ = \ \I_{X/S} \    \phi\cdot
\p_{k <n} c'_1(\cL_k)
\ee
Using induction we get the following  {\it change of metric formula: }
\be
\label{change}
\langle\cL_0\otimes 
\cO(\phi_0),...,\cL_n\otimes \cO(\phi_n)\rangle(\cX/S) \ = \
\langle\cL_0,...,\cL_n\rangle(\cX/S)\otimes \cO(E)
\ee
where
\be
\label{big change formula}  
E \ = \ \I_{\cX/S} \  \s_{j=0}^n  \phi_j\cdot
\p_{k <j} c'_1(\cL_k\otimes \cO(\phi_k))\cdot\p_{k > j} c_1'(\cL_k)
\ee
Let $\phi: \cX\ra \cY$ be a morphism of projective flat integral schemes
over $S$ and let  $m=\dim(\cX/\cY)$ and $n=\dim(\cY/S)$. If
$\cK_0,...,\cK_m$ are line bundles on $ \cX$ and $\cL_1,..,\cL_n$ are line
bundles on $\cY$ then we have a canonical isometry:

$$ \langle
\cK_0,...,\cK_m,\phi^*\cL_1,...,\phi^*\cL_n\rangle(\cX/S)
\ = 
\hskip 1in
$$
\be\label{base change}
\langle\langle
\cK_0,...,\cK_m\rangle(\cX/\cY),\cL_1,...,\cL_n\rangle(\cY/S)
\ee
Now if $\cN$ and $\cN'$ are hermitian line bundles on $S$, then it follows
easily from the definition of the Deligne pairing
that we have the following isometries
\be \label{collapse}
\langle \cL_1,...,\cL_n,\pi^*\cN\rangle(\cY/S) \ = \
\cN^{{\deg}[c_1(\cL_{1\eta})\cdots c_1(\cL_{n\eta})]}
\ee
and 
\be \label{big collapse}
\langle \cL_1,...,\cL_{n-1},\pi^*\cN,\pi^*\cN'\rangle(\cY/S) \ = \ \cO_S
\ee 
Here $L_{j\eta} $ is the restriction of $\cL_j$ to a generic
point of $S$. In particular, if $\cK_0=\phi^*\cL_0$ 
in (\ref{base change}), then the right side reduces to
$\langle \cL_0,...,\cL_n\rangle^{{\deg}[c_1(\cK_{1\eta})\cdots
c_1(\cK_{m\eta})]}$,
where $\cK_{j\eta} $ is the restriction of $\cK_j$ to the fiber at a generic
point of $\cY$.

\subsection{Energy Functionals and the Deligne Pairing}

The change of metric formula (\ref{big change formula}) can be used
to show that various important functionals  have a natural
interpretation as metrics on certain Deligne pairings. This point of view is
rather useful (cf. \cite{PS1} and \cite{PS2} ) and will be adopted as well
in the proofs of the theorems in this paper. In this section we recall
some of the key functionals and give their metric interpretation:
\v
Let $M$ be a K\"ahler manifold with K\"ahler form $\o$. Denote by
$$ 
P(M,\o)\ = \ \{\phi\in C^\i(M,\R): \o_\phi=\o+{\sqrt{-1}\over
2\pi}\pl\bar\pl\phi>0
\
\}
$$
the space of K\"ahler potentials, and
recall that $Ric(\o)=-{\sqrt{-1}\over 2\pi}\ddb\o^n$
is the Ricci curvature form of $\o$. The Aubin-Yau functional is defined by
\be\label{Aubin-Yau functional} {\rm AY}_\o(\phi)\ = \ {1\over [\o^n]}\I_M
\s_{k=0}^n
\phi
\o^{n-k}\o_\phi^k,
\ee
the Futaki functional is defined by
\be \label{Futaki functional}{\rm Fut}_\o(\phi) \ = \ 
\ {1\over [\o^n]}\I_M \s_{k=0}^n \phi
\,Ric(\o)^{n-k}\,Ric({\o_\phi})^k
\ee
and the Mabuchi K-energy is defined by
\bea\label{K-energy}
\n_\o(\phi)\ &=& \ {1\over [\o^n]}
\bigg[
\I_M \left(
\log{\o_\phi^n\over \o^n}
\right)\cdot\o_\phi^n
-\s_{i=0}^{n-1}\I_M\phi Ric(\o)\o_\phi^i\o^{n-1-i} 
\nonumber\\
&&
\quad\quad\quad
+ 
{nc_1(M)[\o]^{n-1}\over (n+1)[\o]^n}\cdot
\s_{i=0}^n \I_M\phi\o_\phi^i\o^{n-i}
\bigg]
\eea
(With a different normalization, the Aubin-Yau functional appears
in the literature as the functional $F_{\o}^0(\phi)=-{1\over n+1}{\rm AY}_{\o}(\phi)$.)
Now suppose that $L\ra M$ is a holomorphic line bundle with metric $h$
and $\o= Ric(h)=-{\sqrt{-1}\over 2\pi}\ddb\,h$ (so in particular, $\o\in c_1(L)$). Corresponding
to the three functions (\ref{Aubin-Yau functional}), (\ref{Futaki functional})
and (\ref{K-energy}), we define the metrized $\Q$ line bundles:
$$ \cA_h\ = \ \langle L,...,L\rangle^{1\over c_1(L)^n}, 
\ \ \ \cF_h = \ \langle
K^{-1}, K^{-1}, ..., K^{-1} \rangle
$$
\be
\label{M}
\cM_h \ = \ \langle K,L,...,L\rangle^{1\over c_1(L)^n}
\langle L,...,L\rangle^{{n\over n+1}{c_1(M)c_1(L)^{n-1}\over
[c_1(L)^n]^2}}
\ee
For example:  $ \cA_{he^{-\phi}}\ = \ \langle
L\otimes O(\phi/2),...,L\otimes O(\phi/2)\rangle^{1\over c_1(L)^n}$.
A direct calculation from the change of metric formula
(\ref{big change formula}) gives 

\v
{\bf Theorem 2.} {\it The Aubin-Yau, Futaki, and Mabuchi
energy functionals arise as changes of metrics
in the corresponding Deligne pairings}
$$ \cA_{he^{-\phi}}\ = \ \cA_h\otimes O(\AY_\o(\phi )/2), \ \ \ 
 \cF_{he^{-\phi}}\ = \ \cF_h\otimes O(\Fut_\o(\phi )/2)
$$
\be
\label{MM} 
\cM_{he^{-\phi}}\ = \ \cM_h\otimes O(\n_\o(\phi)/2)
\ee

\v
\subsection{Proof of Theorem 1.}

For notational simplicity, we give the proof in the case $k=1$.
\v

Fix $h$, a smooth metric on $L=K^{-1}$. If
$\phi$ is a smooth function on $M$, then 
$L(\phi)$
is the bundle $L$  equipped with the metric $he^{-\phi}$.
If $\sigma\in G$, then  $\sigma: L\ra L(\phi_\sigma)$ is an isometry  for
a uniquely defined smooth function $\phi_\sigma$. 
Here ``smooth" means that for every holomorphic map $f$ from a complex
manifold $Y$ to
$M$,
 the pull back metric $f^*h$ is a smooth metric on $f^*L$ and
the pull back function $\phi_\sigma\circ f$ is a smooth function
on $Y$. For example, one can construct such an $h$ via 
formula (\ref{h}).
\v

Now $F(\sigma_t): \langle L,...,L
\rangle^{}
\ra  
\langle L(\phi_t),...,L(\phi_t) \rangle = \langle L,...,L \rangle
(E(t))$ is an isometry, where $\sigma_t=\exp(tX_\R)$. Thus by (\ref{big
change formula}) 
\be \label{E}
 E(t)\ = \ \s_{j=0}^n \I_{M_{reg}} \phi_t \o_t^j\o^{n-j}
\ee
Here $\o_t = \o+{\sqrt{-1}\over 2\pi}\pl\bar\pl \phi_t$ and
$\phi_t=\phi_{\sigma_t}$.
\v

Now let $s\in \langle L,...,L\rangle$. Then  $F(\sigma_t)(s)= \hat
F(\sigma_t)\cdot s$ so $|\hat F(\sigma_t)|\cdot|s|=|F(\sigma_t)(s)| =
|s|\exp(E(t))$. Since $F(\sigma_t)$ is a homomorphism
$R\ra\C^\times$, we have $F(\sigma_t)=\exp(tG(X))$ with
$G(X)\in \C$, so  
$$ E(t)\ = \ \log |F(\sigma_t)| \ = \ tRe(G(X))
$$ 

In other words, $Re(G(X))=E'(t)$. On the other hand, differentiating
both sides of (\ref{E}):
$$ E'(t)\ = \ \s_{j=0}^n \I_{M_{reg}} \dot\phi_t \o_t^j\o^{n-j}\ + \ 
\s_{j=0}^n \I_{M_{reg}} \phi_t {\sqrt{-1}\over
2\pi}\pl\bar\pl \dot\phi_tj\o_t^{j-1}\o^{n-j}
$$
Integrating by parts in the second integral, we obtain
$$ E'(t)\ = \ \I_{M_{reg}}\dot\phi_t \o^n\ = \ Re(F(X))
$$
The last equality is proved using integration by parts (see p. 23 of
\cite{FMS} for the proof for smooth manifolds, and \cite{Y2} for
the general case). We conclude that $Re(F(X))=Re(G(X))$.
Thus $-Im(F(X))=Re(iF(X))=Re(F(iX))=Re(G(iX))$. To prove
Theorem 1, we must show that $G(iX)=iG(X)$, that is, we must
show that $F:G\ra \C^\times$ is a holomorphic function. This
follows from the following:
\v

{\bf Lemma.} {\it The function $F:G\ra \C^\times$ is an algebraic map, that
is, $F$ is locally defined by rational functions.  }
\v

{\it Proof.} We may assume that $G$ is an infinite group.
Let $\langle \ell_0,...,\ell_n\rangle$ be a section of 
$\langle L,...,L\rangle$. Thus $\cap_{j=0}^n div(\ell_j)=\emptyset$.
Let $\sigma\in G$ be in a small neighborhood of ${\bf 1_G}$, the
identity element of $G$. Define
$$ Z_j(\sigma) = \p_{p<j} div(\sigma\ell_p)\p_{p>j}div(\ell_p) = 
\s_k n_{jk}(\sigma)y_{jk}(\sigma) $$ with $n_{jk}(\sigma)\in \Z $
and $y_{jk}(\sigma)$ points in $M$. 
 Then
$$ F(\sigma)^{-1}\ = \ {\langle \sigma\ell_0,...,\sigma\ell_n\rangle
\over
\langle \ell_0,...,\ell_n\rangle}\ = \ \p_{j=0}^n f_j(\sigma)
$$
where 
$$ f_j(\sigma) \ = \ \p_k \left[ {\sigma\ell_j\over
\ell_j}\right](y_{jk}(\sigma))^{n_{jk}(\sigma)}
$$
We claim that each $f_j(\sigma)$ is rational and defined in a
neighborhood of $\sigma = {\bf 1_G}$. To see this, let
$\pi: G\times M\ra M$ be the projection map,
and let $\m: G\times M \ra M$ be the map $(\sigma, m)\ra \sigma(m)$.
Let $Z_j \sub G\times M$ be defined by 
$Z_j =  \p_{p<j} div(\m^*\ell_p)\p_{p>j}div(\pi^*\ell_p) $. Then
$\pi$ restricted to $Z_j$ induces a map $\pi_j: Z_j\ra G$. It
has the property $\pi_j^{-1}(\sigma)=Z_j(\sigma)$. The map $\pi_j$
is generically finite and 
$f_j(\sigma)= {\rm Norm}_{Z_j/G} \left[ {\sigma\ell_j\over
\ell_j}\right]$ is therefore rational, which implies that
$F(\sigma)$ is rational. Since $F$ is a homomorphism, we conclude
that $F$ is everywhere defined, and is therefore regular at
every point $\sigma\in G$. This proves the lemma.

\medskip

Finally, the corollary follows from the theorem of Zhang \cite{Z},
which states that $\langle L^k,...,L^k\rangle = 
\langle L,...,L\rangle^{k^{n+1}}$ is canonically
isomorphic to the Chow line bundle.

\medskip

\section {The Futaki invariant of a complete intersection.}

The Futaki invariant of a complete intersection has
been computed by Lu \cite{Lu} (see as well the paper
of Yotov \cite{Y2}) and is given as follows:
\v

{\bf Theorem 3. (Lu)} {\it
Let $M$ be the $N-s$ dimensional normal Fano variety in $\P^N$
defined by $F_1=\cdots=F_s=0$, where
$F_1,\cdots,F_s$ are homogeneous polynomials of degrees
$d_1,...,d_s$. Let $X=\s_{i,j=0}^N a_{ij}z_i{\pl\over \pl z^j}$ be
a holomorphic vector field on $\P^N$, normalized so that
the matrix $(a_{ij})$ has trace zero. Assume 
$$ XF_i\ = \ \kappa_iF_i, \ \ i=1,...,s
$$
for complex constants $\kappa_1,...,\kappa_s$. Then
$$ F(X)\ = \ (N-s+1)\left[ m^{N-s}\p_{i=1}^s d_i\right]\cdot\left(
-\s_{i=1}^s\kappa_i+{m\over N-s+1}\s_{i=1}^s{\kappa_i\over d_i}
\right)
$$
where $m=N+1-d_1-\cdots -d_s$.
}
\v

We wish to show how Theorem 1 can be used to give
a proof of Theorem~3. 
To do this, it is useful to view $M$ as a fiber in the
family of all complete intersections: 
\v

Let $R_d$ be the space of homogeneous polynomials of degree $d$, and
let
$$ Z \ = \ \P(R_{d_1})\times \cdots \times \P(R_{d_s})
$$
where $s=N-n$. Let
$$ \cX \ = \ \{(z,x): z=(F_1,...,F_s)\in Z, \  x\in \P^N, F_1(x)=\cdots
= F_s(x)=0\}
$$
Then $\cX \sub Z\times \P^N$ is a smooth variety and $dim(\cX)=dim(Z)+n$.
Let $\cL\ra \cX$ be the pull back of the hyperplane bundle on $\P^N$
equipped with the Fubini-Study metric, 
and let $\cK\ra \cX$ be the relative canonical bundle of $\cX/Z$.
Let
$$ 
Z_0\ = \ \{z\in Z: \pi^{-1}z \ \ \ \ {\rm has \ dimension  }\ n\ \}
$$
where $\pi:\cX\ra Z$ is the natural projection. Then $Z_0\sub Z$ is
Zariski open, and the restricted map $\pi:\cX_0\ra Z_0$ is flat
(where $\cX_0=\pi^{-1}(Z_0)$). We want to compute the Deligne pairing
$$
\cF \ = \ \langle \cK^{-1},\cK^{-1},...,\cK^{-1}\rangle(\cX_0/Z_0)
$$
which is a line bundle on $Z_0$. 
\v

Note that $F_j(z)$, when viewed as a function of $(F_j, z)$, is homogeneous of
degree one in $F_j$ and homogeneous of degree $d_j$ in $z$. Thus it represents
a global section $\sigma_j$ of
$H_j\otimes\cL^{d_j}\equiv H_j\cL^{d_j}$. The simultaneous vanishing of the $\sigma_j$ defines the
subvariety $\cX_0\sub Z_0$. Thus we have,
applying (\ref{1.1.5}) successively $s$ times, 
\setcounter{equation}{0}
\be\label{1.1.55} \langle \cK^{-1},...,\cK^{-1}\rangle(\cX_0/Z_0)\ = \
\langle H_1\cL^{d_1},...,H_s\cL^{d_s},
\cK^{-1},...\cK^{-1}\rangle((Z_0\times\P^N)/Z_0)
\ee
where $H_k$ is the pullback of the hyperplane bundle on $\P(R_{d_k})$.
Next we recall that the relative canonical bundle of $Z_0\times \P^N\ra Z_0 $
is $\cL^{N+1}$. Since $\cX$ is a complete intersection, defined by
the sections $\sigma_1\ = \cdots \ = \sigma_d=0$, and since $\sigma_j$
is a global section of $\cL^{d_j}H_j$,
applying the adjunction formula $s$ times yields
$$ \cK^{-1}\ = \ \cL^mH_1^{-1}\cdots H_s^{-1}
$$
where $m=N+1-d_1-\cdots -d_s$.
Expanding  (\ref{1.1.55}) using the multi-linearity of the Deligne
pairing, we get the main term,
$\langle
\cL,...\cL\rangle^{d_1\cdots d_s\cdot m^{n+1}}$, which is the trivial
line bundle: Indeed, $\langle z_0,z_1,...,z_N\rangle $ is a global nowhere
vanishing section of $\langle \cL,...,\cL \rangle$. As
for the other terms in the expansion of (\ref{1.1.55}): If two or more of the
$H_i$ appear, say $H_i$ and $H_j$, we get the trivial bundle, by applying
(\ref{big collapse}), with $\cN=H_i$ and $\cN'=H_j$. Thus the only terms which
contribute are those with one
$H$. There are two ways this can happen: Either the $H_i$ comes from the
$i^{th}$ position in (\ref{1.1.55}), with $1\leq i\leq s$,  and then it is
paired with
$\cL^{d_k}$ with
$k\not=i$, as well as $n+1=N-s+1$ copies of $\cL^m$. Or $H_i$ comes
from any one of the last $n+1$ positions and then it is paired with
$\cL^{d_r}$,
$1\leq r\leq s$ as well as $n=N-s$ copies of $\cL^m$. 
In the first case we apply (\ref{collapse}) with $\cN=H_i$ and
$\{\cL_1,...,\cL_{n}\}\ = \ \{\cL^{d_1},...,\cL^{d_{i-1}}, \cL^{d_{i+1}}, ...,
\cL^{d_{s}},\cL^m, ... , \cL^m
\}$ with $\cL^m$ repeated $N-s+1$ times. 
In the second case we apply (\ref{collapse}) $N-s+1=n+1$ times with
$\cN=H_i^{-1}$ and
$\{\cL_1,...,\cL_{n}\}\ = \ \{\cL^{d_1}, ...,
\cL^{d_{s}},\cL^m, ... , \cL^m
\}$ where $\cL^m$ is repeated $N-s$ times. 
Thus we obtain
\be\label{a_i}
\langle H_1\cL^{d_1},...,H_s\cL^{d_s},
\cL^mH_1^{-1}\cdots H_s^{-1},...,\cL^mH_1^{-1}\cdots H_s^{-1}\rangle\ =
\ \p_{i=1}^s H_i^{a_i}
\ee
where 
\be\label{A} {a_i} \ = \ (-1){(n+1)m^{N-s}d_1d_2\cdots
d_s}\ +\ {m^{N-s+1}}(d_1d_2\cdots d_s)/d_i
\ee

Combining (\ref{1.1.55}) and (\ref{a_i}) we obtain $F(X)=\s a_i\kappa_i$
and this, together with formula (\ref{A}) yields Theorem 3.
\v

\section{Energy functionals for a complete intersection.}

\subsection{The adjunction formula with metrics.}

We wish to prove a non-linear version of Theorem 3 in which
the Mabuchi K-energy is expressed in terms of certain norms on the
polynomials defining $M$. To do this, we first prove a 
metrized version of the adjunction formula:
\v

Let $X$ be a smooth variety and $Y\sub X$
a smooth subvariety. The adjunction formula
says
\setcounter{equation}{0}
\be
\label{5.1} K_Y\ = \ \big(K_X\otimes O(Y)\big)|_Y
\ee
where $O(Y)$ is the line bundle associated
to the divisor $Y$.
\v

The isomorphism is easy to describe in
local coordinates: Let $y\in Y$. Choose
local coordinates $(z_1,...,z_n)$
 in such a way that $y$
is the origin, and $Y$ is the set $z_1=0$.
Then a local section of $K_X\otimes O(Y)$
is meromorphic differential form of the type
$$\eta= {g(z)dz_1\wedge dz_2\wedge \cdots \wedge dz_n\over z_1 }
$$
where $g(z)=g(z_1,...,z_n)$ is a holomorphic function. The isomorphism
(\ref{5.1}) is
$\eta\mapsto \t$ where
$$ \t \ = \ g(0,z_2,...,z_n)dz_2\wedge\cdots dz_n
$$
that is, $\t$ is the {\it residue} of $\eta$. Note that the map $\eta\mapsto \t$
does not depend on the choice of coordinates: If $z_1=uz_1'$ where $u$ is
a nowhere vanishing function, then
${dz_1\over z_1}={dz_1'\over z_1'} + {d(\log u)}$ and thus
$\eta = {g(z')dz_1'\wedge dz_2\wedge \cdots \wedge dz_n\over z_1' } +
{g(z')z_1'd(\log u)\wedge dz_2\wedge \cdots \wedge dz_n\over z_1' }
$.
In other words,
$$ \eta\ = \
{[g(z')+z_1'{\pl \log u\over \pl z_1' }]dz_1'\wedge dz_2\wedge \cdots \wedge dz_n\over z_1' }
$$
where $z'=(z_1',z_2,...,z_n)$.
But $z_1'{\pl \log u\over \pl z_1' }$ vanishes on $Y$. This
shows the isomorphism (\ref{5.1}) is well defined.
\v

An equivalent way of formulating (\ref{5.1})
is
\be \label{5.1'} K_Y\otimes O(-Y)|_Y\ = \ (K_X)|_Y
\ee
Again, this is easy to describe in local
coordinates: Let $y\in Y$, let
$U\sub X$ be an open set containing $y$,
and let $W=U\cap Y \sub Y$. As before, we assume that $Y$ is locally given by
the equation $z_1=0$.  Then a section
of $K_X$ over $U$ is a differential
form $\eta = g(z)dz_1\wedge \cdots
\wedge dz_n$ with $g$ holomorphic on $U$.
We say $\eta$ is equivalent to zero if
$g(0,z_2,...,z_n)=0$.
If $\eta$ and $\eta'$ are sections of
$K_X$ over $U$ and $U'$ respectively,
where $U\cap Y=U'\cap Y = W$, we say that $\eta\sim\eta'$ if
$\eta-\eta'$ is equivalent to zero.
A section of $(K_X)|_Y$ over $W$ is then an equivalence
class $[\eta]$. Similarly, a section
of $O(-Y)$ is a holomorphic function
$f$ on $U$ such that $f(0,z_2,...,z_n)=0$,
in other words, $f$ is divisible by $z_1$.
We say that $f$ is equivalent to zero
if it is divisible by $z_1^2$. A section
of $O(-Y)|_Y$ is then an equivalence class
$[f]$. Finally a section of $K_Y\otimes O(-Y)|_Y$ is locally given by an
expression of the form
$[f]dz_2\wedge\cdots \wedge dz_n$.
And now (\ref{5.1'}) is given by $[\eta]\mapsto [f]dz_2\wedge \cdots \wedge
dz_n$, where $f=z_1g$.
\v

Now we want prove the adjunction formula with
metrics: Let $\o$ be a K\"ahler metric
on $X$. Then $\o^n$ is a metric on $K_X^{-1}$, and $\o^{n-1}$ is a metric on $K_Y$.
We want to construct a  metric
$||\cdot||_{ad}$
on $O(Y)|_Y$ in
such a way that (\ref{5.1}) becomes an isometry.
In other words, for every $y\in Y$, we want $||\cdot||_{ad}$
to have the property
\be
\label{5.2}
\left|\left|{1\over z_1 }\right|\right|_{ad}\ = \ 1
\ee
at the origin, if $(z_1,...,z_n)$ are
local  coordinates on $U\sub X$, centered
at $y\in Y$, satisfying
\v

a) $Y\cap U =\{z: z_1=0\}$.
\vskip .03in
b) $\o={\sqrt{-1}\over 2\pi}\sum_j dz_j\wedge d\bar z_j$
at the origin.
\v

This is equivalent to constructing a
metric $||\cdot||_{ad}$ on $O(-Y)|_Y$ such that
$\left|\left|{z_1 }\right|\right|_{ad}\ = \ 1$
at the origin, if $(z_1,...,z_n)$ satisfy a) and b). We do this
as follows:
\v

Let $t$ be a section of $O(-Y)$ over some
open set $W\sub Y$. Thus
$t$ is  represented by a holomorphic function $T$ on $U\sub X$, 
with $U\cap X=W$, with the property: $T$
vanishes on $Y$.
Here  two holomorphic functions on $U$ are ``equivalent" if their 
difference vanishes to order at least two on $Y$. Define
$$ ||t||_{ad}^2 \ = \
\left({\o^{n-1}\wedge{\sqrt{-1}\over 2\pi}\ddb \left(
{ |T|^2}\right)\over \o^n
 }\right)\bigg|_Y\ = \ ||(\nabla_\o T)|_{_Y}||^2\ = \
||(\pl T)|_{_Y}||^2
$$
where $\nabla_\o T$ is the gradient of
$T$ with respect to the metric $\o$.
This makes sense since $|T|^2$
is a smooth function on $X$.
Thus $||T||_{ad}^2$ is a smooth
function on $Y$.
\v

The following theorem is closely related to Proposition
3.1 of \cite{Lu}, although the statement there looks rather different:
\v

{\bf Theorem 4.}\ {\it  Let $(X,\o)$ be a K\"ahler manifold
and $Y\sub X$ a smooth submanifold of codimension one.
Let $t$ be a section of $O(-Y)|_Y$ over some
open set $W\sub Y$. Thus
$t$ is an equivalence class  represented by a holomorphic function $T$ on
$U\sub X$,
 with $U\cap X=W$, 
and $T|_Y=0$ 
(here  two holomorphic functions on $U$ are ``equivalent" if their 
difference vanishes to order at least two on $Y$). Define
\be \label{5.3}||t||_{ad}^2 \ = \
\left({\o^{n-1}\wedge{\sqrt{-1}\over 2\pi}\ddb \left(
{ |T|^2}\right)\over \o^n
 }\right)\bigg|_Y\ = \ ||(\nabla_\o T)|_{_Y}||^2\ = \
||(\pl T)|_{_Y}||^2
\ee
where $\nabla_\o T$ is the gradient of
$T$ with respect to the K\"ahler form $\o$. Then
$||\cdot ||_{ad}$ is a well defined metric on $O(-Y)|_Y$ and
\be
\label{iso}
K_Y\otimes O(-Y)|_Y\ = \ (K_X)|_Y
\ee
is an isometry with respect to this metric.}
\v

Note that (\ref{5.3}) 
makes sense since $|T|^2$
is a smooth function on $X$.
Thus $||T||_{ad}^2$ is a smooth
function on $Y$.
\v

{\it Proof.} First observe that
the norm $||\cdot||_{ad}$ has a very
simple description in local coordinates:
Let $y\in Y$ and let $(z_1,...,z_n)$ be local coordinates centered at $y$
satisfying a) and b).
We write $T=z_1F$ where $F$ is holomorphic
on $U$. Then
\be 
\label{5.4} 
\left(||t||_{ad}^2\right)(y)\ = \ \left(|F|^2\right)(y)
\ee
Note that (\ref{5.4}) holds only at the point $y$,
whereas (\ref{5.3}) holds on all of $W$.
\v

To establish Theorem 4 we must prove the following:
\v

1) $||t||_{ad}$ depends only on the
equivalence class of $t$.

\vskip .03in

2) If $f$ is a holomorphic function
on $W\sub Y$, then $||ft||_{ad}=|f|\cdot
||t||_{ad}$.

\vskip .03in

3)  $||t||_{ad}\geq 0$ with equality
if and only if $t=0$.

\vskip .03in

4) If $(z_1,...,z_n)$ are  coordinates
centered
at a point $y\in Y$ satisfying a) and b), and if $Y$ is defined
by $z_1=0$, then $||z_1||_{ad}=1$
(or, more precisely, $||t||_{ad}=1$
where $t$ is the section of $O(-Y)$
represented by $z_1$).
\v

Property 1) says
that $||\cdot||_{ad}$ is well defined.
Properties 2) and 3) say that
$||\cdot||_{ad}$ is a genuine metric
on $O(-Y)|_Y$. And property 4) says
that (\ref{iso}) is an isomorphism.
\v

Property 1) can be seen as follows:
Assume $t\sim t'$. Then $t=t'+ z_1^2F$
where $F$ is holomorphic on $U$, and
$Y$ is defined by $z_1=0$. Then
$\pl t = \pl t' + O(z_1)$, so
$\pl t|_Y=\pl t'|_Y$.
Alternatively, property 1) follows
from (\ref{5.4}) since if $t=z_1F$ and $t'=z_1F'$
then $t\sim t'$ if and only if
$F-F'$ is divisible by $z_1$. Since
$z_1$ vanishes at $y$, property 1)
follows.
\v

To prove 2), we note that in calculating
$d (ft)$, the derivatives can land
on $t$ or $f$. But when we restrict to
$Y$, $t$ vanishes. Thus, for the purposes
of calculating $||ft||_{ad}^2$, both
derivatives land on $t$ and $|f|^2$
factors out. Alternatively, property 2)
follows immediately from (\ref{5.4}): If
$t=z_1F$ then $ft=z_1\tilde fF$, where
$\tilde f$ is any holomorphic extension
of $f$ from $W\sub Y$ to $U\sub X$.
\v

To prove 3), we work in local
coordinates: If $t=z_1F$, where $F$
is a holomorphic function on $U$,
$||t||_{ad}=|F(y)|^2$. On the other
hand, $t$ vanishes at $y$ if and
only if $F(y)=0$. This proves 3).
\v

Finally, 4) follows immediately from
(\ref{5.4}), and Theorem 4 is proved.
\v

{\it Remark}: We also have the following relative version of Theorem 4:
Let $S$ be a complex manifold and $\pi: \cX\ra S$  a smooth 
family of complex manifolds of relative dimension $n$. Let $\o$ be a $(1,1)$
form on $\cX$ such that the restriction of $\o$ to each fiber $\cX_s$
is a K\"ahler form. Let ${\cal Y}\sub \cX$ be a smooth submanifold
of codimension one such that $\cY_s={\cal Y}\cap \cX_s$ is smooth
of codimension one. Let $\cK_{\cal Y}$ and $\cK_\cX$ be the relative
canonical bundles of ${\cal Y}\ra S$ and $\cX\ra S$, and endow
$\cK^{-1}_\cX$ with the metric $\o^n$.
Let $t$ be a section of $O(-{\cal Y})|_{\cal Y}$ over some
open set $W\sub\cal Y$. Thus, as before,
$t$ is  represented by a holomorphic function $T$ on $U\sub \cX$, with
$U\cap \cX=W$, with the property: $T$ vanishes on $\cal Y$.
Define $||t||_{ad}^2$ by (\ref{5.3}). Then
$K_{\cal Y}\otimes O(-{\cal Y})|_{\cal Y}\ = \ (K_\cX)|_{\cal Y}
$
is an isometry with respect to this metric. 
\v

\subsection{The Aubin-Yau functional as a norm.}

The Aubin-Yau functional ${\rm AY}_{\o}(\phi)$
on the space $P(M,\o)$ of K\"ahler potentials
was defined in subsection \S 2.3. In this subsection,
we shall show that, for complete intersections, it 
can be written explicitly in terms
of a norm on the space of defining polynomials.
Let $(M,\o)$ be a compact n-dimensional K\"ahler manifold,
and $(L,h)\to M$ a holomorphic line bundle with
metric $h$ satisfying $Ric(h)=\o$. 
The key property which we exploit is the close relation
of ${\rm AY}_{\o}(\phi)$ with the Deligne pairing
\be
{\cal A}_h=\langle L,\cdots,L\rangle^{1\over c_1(L)^n}
\ee
described in Theorem 2, which we reproduce here for convenience
$$ \cA_{he^{-\phi}}\ = \ \cA_h\otimes O({\rm AY}_\o(\phi))
$$
We wish to calculate $\cA_h$ in the case where $M$ is
a complete intersection. First we 
recall our previous notation:
Let $N>n>0$ be positive integers, let $R_d\sub \C[z_0,...,z_N]$ be the space of homogeneous polynomials of degree $d$, and
let
$ Z \ = \ \P(R_{d_1})\times \cdots \times \P(R_{d_s})
$
where $s=N-n$. Let
$$ \cX \ = \ \{(z,x): z=(F_1,...,F_s)\in Z, \  x\in \P^N, F_1(x)=\cdots
= F_s(x)=0\}
$$
Then $\cX \sub Z\times \P^N$ is a smooth variety and $dim(\cX)=dim(Z)+n$.
Let $\cL\ra \cX$ be the pull back of the hyperplane bundle on $\P^N$
equipped with the Fubini-Study metric. Let
$$ Z_0\ = \ \{z\in Z: \pi^{-1}z \ \ \ \ {\rm has \ dimension  }\ n\ \}
$$
where $\pi:\cX\ra Z$ is the natural projection. Then $Z_0\sub Z$ is
Zariski open and the restricted map $\pi:\cX_0\ra Z_0$ is flat
(where $\cX_0=\pi^{-1}(Z_0)$). Then the Deligne pairing
$$
\cS \ = \ \langle \cL,...,\cL\rangle
(\cX_0/\cZ_0)
$$
is a line bundle on $Z_0$ with a continuous hermitian metric.
\v

Let $H$ be the hyperplane bundle on the projective space
$\P(R_{d})$, endowed with the Fubini-Study metric.
A global section of $H$ is
a linear map $R_{d}\ra \C$.
Similarly, a global section of $\cL^d$
is a polynomial map $\C^{N+1}\ra \C$,
homogeneous of degree $d$. And a global
section of $H\otimes \cL^d$ is a
polynomial map $R_d\times \C^{N+1}\ra \C$
which is homogeneous of degree one in the first
variable, and degree $d$ in the second
variable. As before, we observe that
$H\otimes \cL^d\equiv H\cL^d$, which is a line bundle
on $\P(R_d)\times \P^N$, has a canonical global
section $\si_d$ given by the map
$(F,z)\ra f(z)$. The norm of this section
is
$$ ||\si_d||\ = \ {|F(z)|\over |F|\cdot |z|^d }
$$
where $|F(z)|$ and $|z|$ are
the norms of vectors in ${\bf C}^s$ and in ${\bf C}^{n+1}$
respectively, and $|F|$ is the norm of $(F_1,\cdots,F_s)$, viewed
as an element of $R_{d_1}\times\cdots\times R_{d_s}$.
We also observe that
the curvature $c_1'(H\cL^d)$ is
given by $c_1'(H)+dc_1'(\cL)$.
\v

Now we consider the bundle
$$ \G\ = \ \langle H_1\cL^{d_1},...,
H_s\cL^{d_s}, \cL, ..., \cL\rangle
\big((Z_0\times\P^N)/Z_0\big)
$$
We evaluate $\G$ in two different ways:
First we use the multi-linearity of the Deligne pairing
to write $\G$ as a product of various terms: There is one
term that involves none of the $H_k$: As in the proof of
Theorem 2, this term gives
the trivial line bundle over $Z_0$ equipped with a
constant metric.  The terms that involve two $H_k$ are also
trivial, by (\ref{big collapse}) . Thus
$$ \G \ = \ \p_{k=1}^s \langle \cL^{d_1},...,H_k,...,\cL^{d_s}\rangle
$$
But now, applying (\ref{collapse}), we get
$$ \G \ = \ \p_{k=1}^s H_k^{p_k}
$$
where $p_k=\p_{j\not= k}d_j$.
\v

On the other hand, we can evaluate $\G$ in a different way,
applying (\ref{1.2.1}) $s$ times, using $\si_1,...\si_s$ successively:
$$  \G \ = \ \langle \cL,...,\cL\rangle\otimes O(-E)
$$
where $E$ is the function on $Z_0$ defined by
$$  E(F_1,...,F_s)\ = \ \, \s_{k=1}^s \I_{X_{k-1}}
{\log\left({|F_k(z)|\over |F_k|\cdot |z|^{d_k} }\right)\cdot
(d_{k+1}\cdots d_s)\cdot
 c_1'(\cL)^{N-k+1 }
}
$$
and $X_k\ = \{z\in \P^N: F_1(z)=\cdots = F_k(z)=0\}$.
Thus we see that
$\langle \cL,...,\cL \rangle^{-1}=$ $\p_{k=1}^s H_k^{-p_k}\otimes O(-E)$.
Since
$V_k = \I_{X_k} c_1'(\cL)^{N-k} = (d_1\cdots d_k)$ we
conclude that $\langle \cL,...,\cL \rangle^{-1}$ is the
line bundle $\p_{k=1}^s H_k^{-p_k}$ equipped with the
metric
\be\label{aubin} \log||F||_{_{\cal A}}\ = \
\s_{k=1}^s p_k\cdot{1\over V_k }\I_{X_{k-1}}
{\log\left({|F_k(z)|\over |z|^{d_k} }\right)\cdot
c_1'(\cL)^{N-k+1 } }
\ee
where $F=F_1^{\otimes p_1}\otimes \cdots \otimes F_s^{\otimes p_s}$.
Note that if we replace $F_k$ by $\l F_k$, then the right side
is transformed by adding $p_k\log\l$ (which shows that the
right side does indeed define a norm on
$\p_{k=1}^s H_k^{-p_k}$). 

\v

Thus 
$\langle \cL,...,\cL \rangle^{}$ is the
line bundle $\p_{k=1}^s H_k^{p_k}$ equipped with the
metric
\be\label{aubin1} \log||F^*||_{_{\cal A^*}}\ = \
-\s_{k=1}^s p_k\cdot{1\over V_k }\I_{X_{k-1}}
{\log\left({|F_k(z)|\over |z|^{d_k} }\right)\cdot
c_1'(\cL)^{N-k+1 } }
\ee
where $F^*=F_1^{\otimes -p_1}\otimes \cdots \otimes F_s^{\otimes -p_s}$

\v

$\si: \langle \cL\otimes O(\phi_\si/2),...,\cL\otimes
O(\phi_\si/2)\rangle = \langle \cL,...,\cL\rangle \otimes
O(\AY_\o(\phi_\si)/2)
\ra \langle \cL,...,\cL\rangle$ implies   
$ \AY_\o(\phi_\si)/2=  \log||F||_{_{\cal A^*}}- \log||F^\si||_{_{\cal A^*}}
$
$ =  \log||F^\si||_{_{\cal A}}- \log||F||_{_{\cal A}}
$

\v

We now obtain the following:
\v

{\bf Theorem 5. }\ 
{\it Let $M\sub \P^N$ be a complete intersection defined by
the equations $F=0$, where $F=(F_1,...,F_s)$ and
$F_k$ is a homogeneous polynomial of degree $d_k$. Let
$\o$ be the Fubini-Study metric on $M$. Then for all
$\sigma\in GL(N+1)$  we have
\be  \label{aubin1}
{\rm AY}_\o(\phi_\si)\ = \ {1\over [\o]^n}\log{||F^\si||^2_{_{\cal A}}\over
||F||^2_{_{\cal A}} }
\ee
where $\phi_\sigma \ = \ \log{|\sigma x|^2\over |x|^2}$ and
$||\cdot ||_{_{\cal A}} $ is the norm on $\p_{k=1}^s H_k^{-p_k}$
defined by (\ref{aubin}).
}
\v

\subsection{The Mabuchi K-energy as a degenerate norm.}

In this section we show that the Mabuchi $K$-energy of a complete
intersection can be expressed in terms of a certain
degenerate norm on the space of defining polynomials for $M$. 
This may be viewed as a non-linear generalization of Lu's
formula. It also generalizes expressions
for the $K$-energy in \cite{T94} and in \cite{PS1}. 
\v

The Mabuchi $K$-energy was defined in \S 2.3.
With the same notation as in there and in \S 4.2, 
we recall that it is closely related with the 
following Deligne pairing
\be
\cM_h\ = \ \langle K,L,...,L\rangle^{1\over c_1(L)^n}
\langle L,...,L\rangle^{{n\over n+1}{c_1(M)c_1(L)^{n-1}\over
[c_1(L)^n]^2}}
\ee
where $K$ is endowed with the metric $\o^{-n}$.
More precisely, as stated in Theorem 2, we have
\be
\label{changeM} 
\cM_{he^{-\phi}}\ = \ \cM_h\otimes O(\n_\o(\phi))
\ee

Next consider the family
$\cX_0\ra Z_0\sub Z_0\times \P^N$, where $Z_0 =\{F\in Z: F=0$ is a
normal variety of dimension $n\}$. : As before, let
$\cL\ra
\cX_0$ be the pull back of the hyperplane bundle on $\P^N$ to the
manifold
$\cX_0$, and let $\cK^{-1} \ra \X_0$ be the relative anti-canonical
bundle, endowed with the metric $\o_{\P^N}^n$.
Now let
$$ \cM\ = \ \langle \cK,\cL,...,\cL\rangle^{1\over c_1(L)^n}
\langle \cL,...,\cL\rangle^{{n\over n+1}{c_1(M)c_1(L)^{n-1}\over
[c_1(L)^n]^2}}
$$
Let $F\in Z_0$ and let $f\in \cM_F$, where $\cM_F$ is the fiber
of $\cM$ above $F$. The change of metrics formula (\ref{changeM}) implies, for
all $\sigma\in GL(N+1)$
\be
\label{action} 
\n_\o(\phi_\sigma)\ = \ \log{||f^\sigma||_{_\cM}^2\over
||f||_{_\cM}^2}
\ee
 
We wish to calculate
the norm $||\cdot||_{_\cM}$.
The adjunction formula says
\be\label{adj}
\cK^{-1}=\cL^{N+1}O(-Y_1)_{ad}\cdots O(-Y_s)_{ad}
\ee

 where
$Y_i\sub\cX_0$ is the hypersurface defined by $F_i(z)=0$ and
$O(-Y_i)_{ad}$ is the bundle $O(-Y_i)$ endowed with the metric of
Theorem 4. Here we are using the fact that 
\be \label{adj PN}K^{-1}_{\P^N} \ = \ \cO_{\P^N}(N+1)
\ee
Note that the isometry in (\ref{adj PN}) is equivariant
with respect to the action of $SL(N+1)$, and not 
with respect to the action of $GL(N+1)$. It is for this
reason that, unlike the transformation formula (\ref{aubin})
for the Aubin-Yau functional, which holds for all $\sigma\in GL(N+1)$, our
transformation formula for the K-energy will only hold for
$\sigma\in SL(N+1)$.
\v

Since
$F_i(z)$ is a global section of $\cL^{d_i}H_i$, we have a canonical isomorphism
 $(\cL^{d_i}H_i)_{ad} \approx O(-Y_i)_{ad}$ for some uniquely defined
metric $||\cdot ||_{ad} $ on $(\cL^{d_i}H_i)$. Thus
$ \cK^{-1}\ = \ \cL^{N+1}(\cL^{d_1}H_1)_{ad}^{-1}\cdots
(\cL^{d_s}H_s)_{ad}^{-1}
$. 
On the other hand, (\ref{5.3}) implies that
$$ (\cL^dH)^{-1}_{ad}\ = \ (\cL^dH)^{-1}\otimes O
\left(
{ -{1\over 2}
{\log\bigg(
{\o^{n-1}\wedge {\sqrt{-1}\over 2\pi}\ddb{|F(z)|^2\over |z|^{2d}|F|^2}
\over
\o^n
}
\bigg)
}
}
\right)
$$
Thus
\be \label{part}\langle \cK^{-1},\cL,...,\cL\rangle\ = \ 
\langle \cL^{N+1}(\cL^{d_1}H_1)^{-1}\cdots
(\cL^{d_s}H_s)^{-1}, \cL,...,\cL\rangle
\otimes O(-G)
\ee
where $G$ is the function on $\P(R_{d_1})\times \cdots\times \P(R_{d_s})
$ given by
$$ G(F)= G(F_1,..., F_s)\ = \ {1\over 2}\s_{i=1}^s
{ 
\I_{M}{\log\bigg(
{\o^{n-1}\wedge {\sqrt{-1}\over 2\pi}\ddb{|F_i(z)|^2\over |z|^{2d_i}|F_i|^2}
\over
\o^n
}
\bigg)\cdot\o^n
}
}
$$
Here $M$ is the variety $F=0$ and $\o=\o_{FS}$ is the Fubini-Study metric
on
$\P^N$. 
Let 
$$ q_k \ = \ 1-{m\over d_k(n+1)}
$$
where $m=N+1-d_1-\cdots - d_s$.
One sees easily that $q_k\geq 0$ with equality if and only if $d_i=1$ for
all $i$. 
Then (\ref{part}) implies
$$ \cM\ = \ \p_{i=0}^s   H_i^{q_i}
$$
and that $\cM^{-1} $ is equipped with the singular metric
$$ \log||F||_{_\cM}^2\ = \ -{m\over (n+1)d}\cdot\s_{k=1}^s
p_k\cdot{1\over V_k }\I_{X_k} {\log\left({|F_k(z)|^2\over |z|^{2d_k}
}\right)\cdot c_1'(\cL)^{N-k+1 } } \ 
$$
\be\label{sharp}
 +\ {1\over d}\s_{i=1}^s
{ 
\I_{M}{\log\bigg(
{\o^{n-1}\wedge {\sqrt{-1}\over 2\pi}\ddb{|F_i(z)|^2\over |z|^{2d_i}}
\over
\o^n
}
\bigg)\cdot\o^n
}
}
\ee
where $F=(F_1,...,F_s)$. Applying (\ref{action}) we deduce the following
:
\v

{\bf Theorem 6.} \ {\it Let $M\sub\P^N$ be a smooth
complete intersection defined by $F_1=\cdots = F_s=0$, where
$F_j$ is homogeneous of degree $d_j$ . For
$\sigma\in SL(N+1)$ let
$\phi_\si=\log{|\si(z)|^2\over |z|^2}$, with $z\in\P^N$.
Let $d=d_1\cdots d_s$ and $p_k=d/d_k$.
 Let $||F ||_{_\cM}$
be the degenerate metric defined by (\ref{sharp}). Then
for all $\sigma\in SL(N+1)$ we have
$$ \n_\o(\phi_\sigma)\ = \ \log {||F^\sigma||_{_\cM}^2\over ||F||_{_\cM}^2}
$$
}

\v
{\it Remark.} The Aubin-Yau functional and the Mabuchi energy
functional were written respectively in \cite{Z} and \cite{PS1}
in terms of suitable norms of the Chow point of $M$.
When $M$ is a hypersurface, the Chow point of $M$ coincides with
$M$, and the norms $||\cdot||_{_{\cA}}$ and $||\cdot||_{_{\cM}}$
coincide with the norms $||\cdot||$ in \cite{Z}
and $||\cdot||_{\#}$ in \cite{PS1}.
Of course, when $M$ has higher codimension, no such direct
comparison is possible, since the space of Chow points and the
space of defining polynomials are then quite different.
An exact expression for the Mabuchi functional on hypersurfaces
and an asymptotic expression in the general case in terms of Quillen metrics can be found in \cite{T94} and \cite{T97}.
See also \cite{P} for an asymptotic expression for the
Aubin-Yau functional.

\v
{\it Remark.} The approach in the present paper can be extended
to the case of irreducible normal projective varieties.
This together with applications to stability will be reported
elsewhere.

\v\v
{\bf Dedication.} This paper is dedicated to Professor Peter Li,
former editor-in-chief of the Communications in Analysis and Geometry, on the occasion of his fiftieth birthday.

D.H. PHONG 

DEPARTMENT OF MATHEMATICS

COLUMBIA UNIVERSITY, NEW YORK, NY 10027

{\bf phong@math.columbia.edu}

\v

JACOB STURM

DEPARTMENT OF MATHEMATICS

RUTGERS UNIVERSITY, NEWARK, NJ 07102

{\bf sturm@andromeda.rutgers.edu}

\enddocument